\documentclass{amsart}
\usepackage{amsthm,amssymb,amsmath,dsfont,bigdelim}
\usepackage{graphicx}
\usepackage{float}
\usepackage[arc,all]{xy}
\usepackage{longtable}
\usepackage{stackrel}
\newtheorem{thm}{Theorem}[section]
\newtheorem{cor}[thm]{Corollary}

\newtheorem{prop}[thm]{Proposition}
\theoremstyle{definition}
\newtheorem{defn}[thm]{Definition}

\newtheorem*{theorem}{\textbf{Open Problem}}
\newtheorem*{proof M}{\textbf{Proof of  Main Theorem}}
\theoremstyle{remark}

\numberwithin{equation}{section}
\begin{document}
\title[commutativity degree]{On the  values of commutativity degree of    Lie algebras }%
\author{Afsaneh Shamsaki}%
\address{Department of Mathematics, Faculty of Mathematical Sciences,
Ferdowsi University of Mashhad, Mashhad, Iran }
\email{Shamsaki.afsaneh@yahoo.com}%
\address{Department of Mathematics, Faculty of Mathematical Sciences,
Ferdowsi University of Mashhad, Mashhad, Iran }
\author{Ahmad Erfanian }
\email{erfanian@um.ac.ir }
\address{Department of Mathematics, Faculty of Mathematical Sciences,
Ferdowsi University of Mashhad, Mashhad, Iran }
\author{Mohsen Parvizi}
\email{parvizi@um.ac.ir }
%\thanks{0}%
%\keywords{Commutativity degree, Asymptotic commutativity degree, Lie algebra}%
%\thanks{0}%
%%\thanks{0}%
%\keywords{exterior degree, }%
%\keywords{Schur multiplier, nilpotent Lie algebra}

%\date{0}%
%\dedicatory{0}%
%\commby{0}%
% ----------------------------------------------------------------
\begin{abstract}
In this paper, the possible values of commutativity degree of Lie algebras are determined. Also,
we define the asymptotic commutativity degree of  Lie algebras and obtain the asymptotic commutativity degree for  some of them. Moreover, we prove the existence of a family of Lie algebras such that   the asymptotic commutativity degree is equal to $ \frac{1}{q^{k}} $ for all $ q\geq 2 $ and  a positive integer $ k. $
\end{abstract}
\maketitle
\section{ Introduction}
In the last few years, combining branches of mathematics have attracted  some authors. One of the exciting research topics is the relation between group theory and probability theory. The probability that two elements of a group commute is the commutativity degree of a group. Several authors have  published about it and some interesting  facts have been obtained.
 For instance, it was proved  that there is no group with the commutativity  degree in the interval $ (\frac{5}{8}, 1). $  Lescot in \cite{Lescot} showed that if $ G/Z(G)\cong S_3, $ 
then $ d(G)=\frac{1}{2}. $  An important concept about the commutativity degree of groups is called the asymptotic commutativity degree
of groups. Suppose that $ G_n $ is a family of finite non-abelian groups such that $ | G_n|\rightarrow \infty $ as $ n\rightarrow \infty. $ Then the limit of $ d(G_n) $ as $ n\rightarrow \infty  $ is the asymptotic commutativity degree of $ G_n. $ In  2019, it was shown 
%%Rajat?? \cite{•} 
that there is  a family of finite non-abelian groups whose asymptotic commutativity degree is  $ \frac{1}{k}$ for every positive integer $ k. $ Since group theory and Lie theory have a close relation to each other,  analogous to groups the commutativity degree of finite-dimensional Lie algebra $ L $ in \cite{Sh1} was defined  as follows:
\begin{equation*}
d(L)=\dfrac{| \lbrace (x, y)\in L\times L \mid [x, y]=0\rbrace |}{| L|^{2}}.
\end{equation*}
It is obvious that $d(L)=1  $ if and only if $ L $ is an abelian Lie algebra. 
It can be  shown that $ d(L)=\frac{1}{| L|^{2}}\sum_{x\in L} | C_{L}(x)|,$ where  $ C_{L}(x)=\lbrace y\in L \mid [x, y]=0, ~~ y\in L \rbrace.$ 
Also, we may rewrite $ d(L) $ as the formula 
$ d(L)=\frac{1}{| L|}\sum_{x\in L} \frac{1}{|  \operatorname{Im}ad_{x}|}$ in which  $ ad_x : L\rightarrow L $ is the linear map given by $ y\mapsto [x, y].$ 
   \\
 Determining the possible values of commutativity degree of a Lie algebra is an interesting question. The authors
 determined
 a few values $ d(L) $ in \cite{Sh1}. 
 Moreover, they showed that there is no
   finite-dimensional Lie algebra over the field $ \mathbb{F}_{q} $ with the commutativity degree in the interval $ (\frac{5}{8}, 1). $ Also,
all Lie algebras with $ d(L)= \frac{q^{2}+q-1}{q^{3}}$  have been classified and    upper and lower  bounds have been given  for $ d(L). $
  \\
In this paper, we focus on  possible  values of the commutativity degree of Lie algebras  
  in the interval $ [\frac{1}{q}, \frac{q^{2}+q-1}{q^{3}}].$   
In addition, we  define the asymptotic commutativity degree of Lie algebras and prove the existence of a family of finite-dimensional non-abelian
Lie algebras with the asymptotic commutativity degree $ \frac{1}{q^{k}} $ for every positive integer $ k. $
%show  the only Lie algebra with the derived subalgebra of dimension one has such values.  
\\
Throughout this paper, we always assume that $ L $ is a non-abelian  finite-dimensional Lie algebra over the field $ \mathbb{F}_{q} $ for all $ q\geq 2. $ In the next section, we use the Lie algebras $ A(n) $ and $ H(m)=\langle x_i, y_i, z\mid [x_i, y_i]=z, 1\leq i \leq m \rangle $   which denote an abelian Lie algebra of dimension $ n $ and  Heisenberg  Lie algebra of dimension $ 2m+1, $ respectively. 
%The later proposition is an important equipment in the next section.
%\begin{prop}\label{direct}\cite[Proposition 3.5]{Sh1}
%Suppose that $ H $ and $ K $ are Lie algebras and $ L=H \oplus K. $ Then
%$
%P(L)=P(H) P(K).
%$
%\end{prop}
\section{the Possible values of  commutativity degree of Lie algebras}
In this section, we are going to determine  the possible values of commutativity degree of Lie algebras  in the interval $ [\frac{1}{q}, \frac{q^{2}+q-1}{q^{3}}] $  and obtain the structure of such Lie algebras.
\\
The  following proposition holds for Lie algebra over an arbitrary field.
\begin{prop}\label{strac}
Let $ L $ be an $ n $-dimensional Lie algebra with the derived subalgebra of dimension $ 1. $ Then $ L $ is isomorphic to one of the Lie algebras $ H(m)\oplus A(n-2m-1) $ for all $ m\geq 1 $ or $ \langle x, y\mid [x, y]=x \rangle \oplus A(n-2). $
\end{prop}
\begin{proof}
We consider  the following  cases:\\
\textbf{Case 1.}  $ L^{2}\subseteq Z(L). $ In this case,  $ L $ is of nilpotency  class two. Moreover, $ \dim L^{2}=1 $ implies that  $ L $ is isomorphic to 
$ H(m)\oplus A(n-2m-1) $ for all $ m\geq 1 $ by \cite[Lemma 3.3]{N}.
\\
\textbf{Case 2.}  $ L^{2}\nsubseteq Z(L). $ By taking the non-zero element $ x\in L^{2}\setminus Z(L), $ 
   there is an element $ y\in L $ such that $ [x, y]\neq 0. $
%   Put the set  $ \lbrace x \rbrace $ as a basis for $ L^{2} $ by the assumption $ \dim L^{2}=1. $ 
Since $ \dim L^{2}=1 $ and $ x $ is non-zero,  the singleton set  $ \lbrace x \rbrace $ is in fact   a basis for $ L^{2}. $ 
Hence $ [x, y]=\alpha x $ for some $ \alpha \in \mathbb{F}. $ Without loss of generality, we may assume that
%By replacing $ y $ with a scalar multiplication of itself, we may assume
$ [x, y]=x. $  On the other hand,   the set $ \lbrace x, y \rbrace $ is linearly independent due to  $ [x, y] \neq 0. $ 
 Now, extend $ \lbrace x, y \rbrace $ to a basis $ \lbrace x, y, y_j \mid 1\leq j \leq n-2 \rbrace $ of $ L. $ Since $ x $ generates $ L^{2}, $ there are  scalars $ a_i, b_i, c_{ij} \in \mathbb{F} $ such that 
\begin{equation*}
 [x,y]=x,~~ [x, y_i]=a_i x,~~ [y, y_i]=b_i x,~~ [y_i, y_j]=c_{ij} x
\end{equation*}
 for all $ 1\leq i,j \leq n-2. $ 
 If $ c_{ij}\neq 0 $  for some $ 1\leq i,j \leq n-2, $ then 
\begin{equation}\label{e5}
x=[x,y]=\frac{1}{c_{ij}}[[y_i, y_j],y]=- \frac{1}{c_{ij}}([[y, y_i], y_j]+[[y_j, y], y_i])=- \frac{1}{c_{ij}}(b_ia_j-b_ja_i)x
\end{equation}
%for all  $ 1\leq i< j \leq n-2 $ 
by  the Jacobi identity. Hence $ c_{ij}= b_ja_i -b_ia_j. $
% for all  $ 1\leq i< j \leq n-2. $
Now if $ c_{ij}= 0 $ for some $ i, j $ such that $ 1\leq i,j \leq n-2, $ then  similarity
\begin{equation}\label{e6}
0=[[y_i, y_j],y]=- ([[y, y_i], y_j]+[[y_j, y], y_i])=(b_ja_i -b_ia_j)x
\end{equation}
for all  $ 1\leq i< j \leq n-2 $ and so $  b_ia_j=b_ja_i. $
We claim that the set $ \lbrace z_{t}=-b_tx+a_ty-y_t \mid 1\leq t \leq n-2 \rbrace $ is central. It is clear that 
\begin{align*}
[x, z_t]=0 \quad \text{and} \quad [y, z_t]=0 \quad \text{for~~ all} ~~ 1\leq t \leq n-2.
\end{align*}
Also, $ [y_j, z_t ]=(b_t a_j-b_j a_t)x-[y_j, y_t]=(b_t a_j- b_ja_t-c_{jt})x =0$ for $ 1\leq t\leq n-2 $ by  \eqref{e5} and \eqref{e6}. 
%for all $ a_i, b_i, c_{ij} \in \mathbb{F}_q $ such that $ 1\leq i< j \leq n-2 $ and $ 1\leq  t\leq n-2. $ 
 It can be shown that
 the  set $ \lbrace  x, y, z_t\mid 1\leq t\leq n-2\rbrace $ is linearly independent so, it is  a basis of $ L. $  Hence 
$ L=\langle  x, y, z_t\mid [x, y]=x, ~~ 1\leq t\leq n-2\rangle $ which is
%\begin{equation*}
%L=\langle  x, y, z_t\mid [x, y]=x, ~~ 1\leq t\leq n-2\rangle.
%\end{equation*}
 $ L $ is isomorphic to  $ \langle x, y\mid [x, y]=x \rangle \oplus A(n-2). $
\end{proof}
In the next proposition, we observe that the possible values of commutativity degree of Lie algebras with the derived subalgebra of dimension at least two are strictly smaller than $\frac{1}{q}$ for all $ q\geq 2. $
\begin{thm}\label{2}
Let $ L $ be a Lie algebra with the derived subalgebra of dimension at least two. Then $ d(L)< \frac{1}{q} $ for all $ q\geq 2. $
\end{thm}
\begin{proof}
Since  $ \dim L^{2}\geq 2,$ there is an element $ x\in L $ such that $   \dim \operatorname{Im}ad_{x}=k\geq 2.$ Otherwise,
if $  \dim \operatorname{Im}ad_{x}=1 $ for all $ x\in L,$  then  $ \dim L^{2}=1 $ by \cite[Theorem 2.3]{breath}, which is a contradiction. 
Hence the set $ \lbrace [x, x_i] \mid x_i \in L, 1\leq i \leq k \rbrace  $ is linearly independent. 
Now, we claim that $ d(L) < \frac{1}{q}$ by considering the following cases.\\
\textbf{Case 1.} Suppose that there is an element $ x_i $ such that   $ \dim \operatorname{Im}ad_{x_i}=k_1 \geq 2$ and $1\leq i \leq k. $ 
Then  $  \dim \operatorname{Im}ad_{\gamma_1 x_i+z}=k_1\geq 2  $ for all $z\in Z(L) $ and $ \gamma_1 \in \mathbb{F}^{*}_q.$ So,  
the number of these elements are  $ (q-1)| Z(L)|. $ If there exists other element $ x_j $ such that 
 $ \dim \operatorname{Im}ad_{x_j}=k_2 \geq 2$ and $1\leq  i\neq j \leq k, $ then  
$  \dim \operatorname{Im}ad_{\gamma_2 x_j+z}=k_2\geq 2  $ for all $z\in Z(L) $ and  $ \gamma_2 \in \mathbb{F}^{*}_q$ and the number of them is $ (q-1) | Z(L)|.$ Also,
$  \dim \operatorname{Im}ad_{\gamma x+z}=k \geq 2  $ for all $z\in Z(L) $ and $ \gamma \in \mathbb{F}^{*}_q$ and the number of
these elements is $ (q-1)| Z(L)|. $ 
%Therefore the number of these elements are  $ (q-1)| Z(L)|. $
 Hence
\begin{align*}
&d(L)=\dfrac{| Z(L)|}{| L|}+
\dfrac{1}{| L|} \sum_{x\in L\setminus  Z(L)} \frac{1}{| \operatorname{Im}ad_{x}|} \cr
&\leq \dfrac{| Z(L)|}{| L|}+\dfrac{(q-1)| Z(L)|}{q^{k}| L|}+\dfrac{(q-1)| Z(L)|}{q^{k_1}| L|}+\dfrac{(q-1)| Z(L)|}{q^{k_2}| L|}\cr
&+\dfrac{| L|-| Z(L)|-3(q-1)| Z(L)|}{q| L|}\cr
&\leq \dfrac{| Z(L)|}{| L|}+\dfrac{(3(q-1)-q-3q(q-1)) | Z(L)|}{| L|q^{2}}+\dfrac{1}{q}\cr
&=\dfrac{| Z(L)|}{| L|} (\dfrac{-2q^{2}+5q-3}{q^{2}})+\dfrac{1}{q}.
\end{align*}
Since $ -2q^{2}+5q-3< 0 $ for all $ q\geq 2, $ we have $ d(L)<\dfrac{1}{q}.  $\\
If  $ \dim \operatorname{Im}ad_{x_j}=1 $ for all $1\leq  i\neq j \leq k, $ then 
$ \dim \operatorname{Im}ad_{\gamma_1 x+\gamma_2 x_j+z}=k_3 \geq 2 $ for all $1\leq  i, j\leq k$  such that $ \gamma_1, \gamma_2 \in \mathbb{F}^{*}_q$ and the number of them is $ (q-1)^{2}| Z(L)| $ when 
 $ [x_i, x_j]=0. $ Therefore
 \begin{align*}
&d(L)=\dfrac{| Z(L)|}{| L|}+
\dfrac{1}{| L|} \sum_{x\in L\setminus  Z(L)} \frac{1}{| \operatorname{Im}ad_{x}|} \cr
&\leq \dfrac{| Z(L)|}{| L|}+\dfrac{(q-1)| Z(L)|}{q^{k}| L|}+\dfrac{(q-1)| Z(L)|}{q^{k_1}| L|}+\dfrac{(q-1)^{2}| Z(L)|}{q^{k_3}| L|}\cr
&+\dfrac{| L|-| Z(L)|-(2(q-1)+(q-1)^{2})| Z(L)|}{q| L|}\cr
&\leq \dfrac{| Z(L)|}{| L|} (\dfrac{-q^{3}+2q^{2}-1}{q^{2}})+\dfrac{1}{q}<  \dfrac{1}{q}.
\end{align*}
Also, we have $ \dim \operatorname{Im}ad_{\gamma_1 x_i+\gamma_2 x_j+z}=k_4 \geq 2 $
for all $1\leq  i\neq j\leq k$  such that $ \gamma_1, \gamma_2 \in \mathbb{F}^{*}_q$ with the number of $ (q-1)^{2}| Z(L)| $ when $ [x_i, x_j]\neq 0. $ 
By a similar way, one can check $ d(L)<\frac{1}{q}. $
\\
\textbf{Case 2.} Assume that $ \dim \operatorname{Im}ad_{x_i}=1 $ for all $1\leq  i\leq k. $  It is easy to see that $ [x_i, x_j]=0 $
for all $1\leq  i, j \leq k. $ In this case $ \dim \operatorname{Im}ad_{\gamma_1 x+z}=k \geq 2 $ and 
$ \dim \operatorname{Im}ad_{\gamma_1 x+\gamma_2 x_i+z}=k_5 \geq 2 $ for all $1\leq  i\leq k$  such that $ \gamma_1, \gamma_2 \in \mathbb{F}^{*}_q.$
So, the number of this elements is $ (q-1)| Z(L)| $ and  $ k(q-1)^{2}| Z(L)|, $  respectively. Hence 
 \begin{align*}
&d(L)=\dfrac{| Z(L)|}{| L|}+
\dfrac{1}{| L|} \sum_{x\in L\setminus  Z(L)} \frac{1}{| \operatorname{Im}ad_{x}|} \cr
&\leq \dfrac{| Z(L)|}{| L|}+\dfrac{(q-1)| Z(L)|}{q^{k}| L|}+\dfrac{k(q-1)^{2}| Z(L)|}{q^{k_5}| L|}
+\dfrac{|L|-((q-1)+1+k(q-1)^{2}) | Z(L)|}{| L|q}\cr
&\leq \dfrac{| Z(L)|}{| L|}+\dfrac{(q-1)| Z(L)|}{q^{2}| L|}+\dfrac{k(q-1)^{2}| Z(L)|}{q^{2}| L|}
+\dfrac{|L|-((q-1)+1+k(q-1)^{2}) | Z(L)|}{| L|q}\cr
&=\dfrac{| Z(L)|}{| L|} (\dfrac{q-1}{q^{2}})(1-k(q-1)^{2})+\dfrac{1}{q}< \dfrac{1}{q}.
\end{align*}
\end{proof}
In \cite{Les}, Lescot proved that all groups with $ d(G)> \frac{1}{2} $ are nilpotent. 
In the later theorem, we prove that 
there is a non-nilpotent  Lie algebra $L$ with  $d(L)> \frac{1}{q},$  which shows  that the above fact does not hold
for Lie algebras. 
\begin{thm}\label{1}
Let $ L $ be an $ n $-dimensional Lie algebra. Then $ d(L)> \frac{1}{q} $ if and only if $ L $ is isomorphic to one of the Lie algebras 
$ H(m)\oplus A(n-2m-1) $ such that $ m\geq 1 $ or $ \langle x, y \mid [x, y]=x \rangle \oplus A(n-2). $
\end{thm}
\begin{proof}
Let $ d(L)> \frac{1}{q}. $  If   $ \dim L^{2}\geq 2, $ then $ d(L)< \frac{1}{q} $ by  Proposition \ref{2}, which is a contradiction. Thus  $ \dim L^{2}=1. $  Hence $ L $ is isomorphic to 
$ H(m)\oplus A(n-2m-1) $ for  $ m\geq 1 $ or $ \langle x, y \mid [x, y]=x \rangle \oplus A(n-2) $ by Proposition \ref{strac} and we know that  the commutativity of these structures is strictly bigger than  $  \frac{1}{q} $  by \cite[Proposition 3.5, Examples 3.3 and 3.4]{Sh1}.   
\\
The converse of theorem is obtained by \cite[Proposition 3.5, Examples 3.3 and 3.4]{Sh1}.   
\end{proof}
%The next corollaries are the results of the above theorem.\\
From \cite[Theorem 5]{Nath}, we know that  for every positive integer $ n $ there is a finite group $ G $ with $ d(G)=\frac{1}{n}. $  By Proposition \ref{2} and Theorem \ref{1}, one can see that there is no Lie algebra over the field $ \mathbb{F}_{q} $ with $ d(L)=\frac{1}{q}. $
Furthermore,
Proposition \ref{strac}, Theorem \ref{1} and  \cite[Corollary 4.3]{Sh1} imply that  $\dim L^{2}=1$ if and only if   $d(L)\in (\frac{1}{q}, \frac{q^{2}+q-1}{q^{3}}]. $
%\\
%The later corollary illustrates that only values of $ d(L) $ in the interval $ [ \frac{1}{q}, \frac{q^{2}+q-1}{q^{3}}] $  are exactly $ \dots, \frac{q^{2n}+q-1}{q^{2n+1}}, \dots,  \frac{q^{8}+q-1}{q^{9}},  \frac{q^{6}+q-1}{q^{7}},  \frac{q^{4}+q-1}{q^{5}}, \frac{q^{2}+q-1}{q^{3}}.$ 
%for all $ q\geq 2. $
\\Now, it is interesting to ask which rational numbers in the interval $[\frac{1}{q}, \frac{q^{2}+q-1}{q^{3}}]$ can be 
the commutativity degree of  Lie algebras with the derived subalgebra of dimension $ 1. $ In the next corollary, we give a complete answer
to this question.
\begin{cor}
Let $ L $ be a Lie algebra such that $ d(L)\in  [\frac{1}{q}, \frac{q^{2}+q-1}{q^{3}}].$ Then  $ d(L) $ is  one of  the  terms of sequence  $\lbrace \frac{q^{2n}+q-1}{q^{2n+1}} \rbrace_{n=1}^{\infty}. $
% for all $ q\geq 2. $
\end{cor}
\begin{proof}
The proof is obtained by Theorem \ref{1}, \cite[ Examples 3.3 and 3.4 ]{Sh1}.
%By Theorem \ref{1}, the Lie algebras  $ H(m)\oplus A(n-2m-1) $ for all $ m\geq 1 $ and $ \langle x, y \mid [x, y]=x \rangle \oplus A(n-2) $ are only structures with $ d(L) \in [\frac{1}{q}, \frac{q^{2}+q-1}{q^{3}}]. $ 
%for all $ q\geq 2. $ 
%On the other hand, $ d(L) $ is equal to one of the elements in the sequence  $\lbrace \frac{q^{2n}+q-1}{q^{2n+1}} \rbrace_{n=1}^{\infty} $ %for all $ q\geq 2 $
%by  \cite[ Examples 2.3 and 2.4 ]{Sh1}.
\end{proof}
\begin{cor}\label{aso1}
Let $L$ be a Lie algebra. Then $$\lbrace d(L) \mid \dim L^{2}=1 \rbrace=\lbrace \frac{q^{2n}+q-1}{q^{2n+1}} \mid n \in \mathbb{N} \rbrace.$$  
\end{cor}
%\begin{proof}
%Let $ L $ be an arbitrary $ n $-dimensional Lie algebra with  $ \dim L^{2}=1.$ Then $ L $ would be isomorphic to one of the Lie algebras $ H(m)\oplus A(n-2m-1) $ for all $ m\geq 1 $ or $ \langle x, y \mid [x, y]=x \rangle \oplus A(n-2), $ by Proposition \ref{strac}. We know that  $ d(H(m)\oplus A(n-2m-1) )= \frac{q^{2m}+q-1}{q^{2m+1}}$ for all $ m\geq 1 $ and $ d( \langle x, y \mid [x, y]=x \rangle \oplus A(n-2))=\frac{q^{2}+q-1}{q^{3}}, $
%by  \cite[Proposition 3.5, Examples 2.3 and 2.4]{Sh1}.    So, $ d(L)$ belongs to the set of right hand side of the above equality. Conversly, if  $L$ is a Lie algebra with $ d(L)=\frac{q^{2m}+q-1}{q^{2m+1}} $ for some $m\geq 1.$
%then we may put  $ L=H(m). $ Therefore $ d(H(m))=\frac{q^{2m}+q-1}{q^{2m+1}} $   and  $\dim H(m)^{2}=1$
%by Proposition \ref{strac}. Hence the proof is completed.
%for every number $\frac{q^{2m}+q-1}{q^{2m+1}}  $  such that $ m\in \mathbb{N}. $ ???
%\end{proof}
%If $ L $ is a Lie algebra over the field $ \mathbb{F}_2, $ then only values of $ d(L) $ in the interval $ [ \frac{1}{2}, \frac{5}{8}] $ are precisely 
%$ \dots,  \frac{2^{2n}+1}{2^{2n+1}}, \dots, \frac{257}{2^{9}}, \frac{65}{2^{7}}, \frac{17}{2^{5}}, \frac{5}{2^{3}}. $
\begin{cor}
Let $ L $ be a Lie algebra over the field $ \mathbb{F}_{2} $ and $ d(L)\in  [\frac{1}{2}, \frac{5}{8}].$ Then $ d(L) $ is one of the terms of  the sequence $ \lbrace \frac{2^{2n}+1}{2^{2n+1}}\rbrace_{n=1}^{\infty}. $
\end{cor}
In \cite[Corollary 4.6]{Sh1}, we determined $ d(L) $ when $ \dim L/Z(L)=2. $ The next theorem  demonstrates that all possible values of $ d(L) $ are two fixed values   provided that  $ \dim L/Z(L)=3. $ 
\begin{thm}\label{central}
Let $ L $ be a Lie algebra such that $ \dim L/Z(L)=3. $ Then $ d(L) $ is  one of the values $ \frac{2q^{2}-1}{q^{4}} $ or $ \frac{q^{3}+q^{2}-1}{q^{5}}. $ 
\end{thm}
\begin{proof}
Assume that $ \dim L/Z(L)=3. $ Then $ L=\langle x_1, x_2, x_3, Z(L) \rangle$ and so $ L^{2}=\langle [x_1, x_2], [x_1, x_3], [x_2, x_3] \rangle. $
Hence $ \dim L^{2}\leq 3.$ On the other hand, $ L $ is non-abelian thus $ \dim L^{2}\neq 0. $ If $ \dim L^{2}=1, $ then $ \dim L/Z(L) $ is
 even by the structures of Lie algebras in Proposition \ref{strac}. It is a contradiction with $ \dim L/Z(L)=3. $ 
 Hence $ \dim L^{2} $ is  $ 2 $ or $ 3. $ Let $ \dim Z(L)=t. $
We have  the following two cases.\\
\textbf{Case 1.}  $ \dim L^{2}=2. $ \\
Without loss of generality, let $ \lbrace [x_1, x_2], [x_1, x_3] \rbrace $ be a basis of $ L^{2}$  and $[x_2, x_3]=\gamma_1 [x_1, x_2]+ \gamma_2 [x_1, x_3]  $ such that $\gamma_1, \gamma_2\in \mathbb{F}_q.$ Also,
 $ a=\alpha_1 x_1+\alpha_2 x_2+\alpha_3 x_3 $  for all $ a\in L, $  $ z \in Z(L) $ and $ \alpha_i \in \mathbb{F}_{q} $ such that  $ 1\leq i \leq 3. $  
 Then
 \begin{align*}
&[a, x_1]=-\alpha_2  [x_1, x_2]-\alpha_3  [x_1, x_3],\cr 
&[a, x_2]=(\alpha_1 -\alpha_3 \gamma_1) [x_1, x_2]-\alpha_3 \gamma_2 [x_1, x_3],\cr 
& [a, x_3]=\alpha_2 \gamma_1 [x_1, x_2]+(\alpha_1+\alpha_2 \gamma_2) [x_1, x_3]. 
\end{align*}
  Since $ \operatorname{Im}ad_{a}\subseteq L^{2} $ for all $ a\in L,  $ we have $ \dim \operatorname{Im}ad_{a} \leq 2. $ \\
\textbf{Subcase 1-1.}   $ \gamma_1=\gamma_2=0. $\\
  If $ \alpha_1\neq 0, $ then $ [a, x_2]=\alpha_1 [x_1, x_2] $ and $ [a, x_3]=\alpha_1 [x_1, x_3]. $ By assumption the set $ \lbrace [x_1, x_2], [x_1, x_3] \rbrace $ is linearly independent, hence $ \dim \operatorname{Im}ad_{a} =2 $ when $ \alpha_1 \neq 0. $
  So, the number of the elements is equal to  $ (q-1)q^{2+t}. $ Otherwise,  $ \dim \operatorname{Im}ad_{a}=1  $ and the number of them is $ (q^{2}-1)q^{t}. $
 %By a similar method, if $ a=\alpha_2 x_2+\alpha_3 x_3+\beta z  $ and $ \alpha_2\neq 0 $ or 
%$ a=\alpha_3 x_3+\beta z  $ such that $ \alpha_3 \neq 0 $  and $ \alpha_2, \alpha_3 \in \mathbb{F}_{q}, $
%then    for the both of cases.
%Hence  the number of such elements are equal to $ (q-1)q^{1+t} $ or $ (q-1)q^{t}, $ respectively. 
Therefore 
\begin{align*}
d(L)&=\dfrac{| Z(L)|}{| L|}+
\dfrac{1}{| L|} \sum_{x\in L\setminus  Z(L)} \frac{1}{| \operatorname{Im}ad_{x}|} \cr
&=\dfrac{q^{t}}{q^{3+t}}+\dfrac{1}{q^{3+t}}\Big(\dfrac{(q-1)q^{2+t}}{q^{2}}+\dfrac{(q^{2}-1) q^{t}}{q}\Big)=\dfrac{2q^{2}-1}{q^{4}}.
\end{align*}
\textbf{Subcase 1-2.} 
 $   \gamma_1 \neq 0$ and $\gamma_2 \neq 0.$ \\ If
 $ \alpha_1 -\alpha_3 \gamma_1=0 $ and $ \alpha_1+\alpha_2 \gamma_2=0,$ then $ a=\alpha_1 x_1-\gamma_2^{-1}\alpha_1 x_2+\gamma_1^{-1}\alpha_1 x_3+z. $ On the other hand, $ a \neq 0 $ thus $ \alpha_1 \neq 0. $ So, $ \alpha_2 $ and $ \alpha_3 $ are
non-zero. Therefore
$ \operatorname{Im}ad_{a}=\langle \alpha_3 \gamma_2 [x_1, x_3],  \gamma_1 \alpha_2  [x_1, x_2] \rangle $
and  $ \dim \operatorname{Im}ad_{a} =2. $
Let $ \alpha_1 -\alpha_3 \gamma_1\neq 0 $ and $ \alpha_1+\alpha_2 \gamma_2= 0.$ Then
$ \alpha_1 \neq \alpha_3 \gamma_1 $ and $ \alpha_1=-\alpha_2 \gamma_2. $ Thus 
$ a=\alpha_1 x_1-\alpha_1 \gamma_2^{-1} x_2+\alpha_3 x_3+z. $ Since $ a\in L\setminus Z(L), $ then $ \alpha_1 $ and $ \alpha_3 $ can not be 
zero together. Therefore  $ \operatorname{Im}ad_{a}=1$ when $ \alpha_1 \neq 0 $ and $ \alpha_3=0. $ Otherwise,  $ \operatorname{Im}ad_{a}=2.$ Suppose that $ \alpha_1 -\alpha_3 \gamma_1=0 $ and $ \alpha_1+\alpha_2 \gamma_2\neq 0.$
By a similar way of previous case, if $ \alpha_1\neq 0 $ and $ \alpha_2=0, $ then $ \operatorname{Im}ad_{a}=1$ and otherwise,
$ \operatorname{Im}ad_{a}=2.$ 
Finally, let $ \alpha_1 -\alpha_3 \gamma_1\neq 0 $ and $ \alpha_1+\alpha_2 \gamma_2\neq 0.$ Then 
 $ \operatorname{Im}ad_{a}=2$ when $ \alpha_2=0 $ or $ \alpha_3=0. $ Assume that   $ \alpha_2 $  and  $ \alpha_3 $ are non-zero.
One can see that $ \operatorname{Im}ad_{a}=1$ when $ \alpha_1=-\alpha_2 \gamma_2+\alpha_3 \gamma_1 $ and if
 $ \alpha_1\neq -\alpha_2 \gamma_2+\alpha_3 \gamma_1, $ then $ \operatorname{Im}ad_{a}=2.$
Now,  the number of elements such $ \operatorname{Im}ad_{a}=1$ is equal to $ q^{2+t}-q^{t}. $ 
Hence 
\begin{align*}
d(L)&=\dfrac{| Z(L)|}{| L|}+
\dfrac{1}{| L|} \sum_{x\in L\setminus  Z(L)} \frac{1}{| \operatorname{Im}ad_{x}|} \cr
&=\dfrac{q^{t}}{q^{3+t}}+\dfrac{1}{q^{3+t}}(\dfrac{q^{2+t}-q^{t}}{q}+\dfrac{q^{3+t}-q^{2+t}}{q^{2}})=\dfrac{2q^{2}-1}{q^{4}}.
\end{align*}
\textbf{Subcase 1-3.} 
By a similar way, $ d(L)=\dfrac{2q^{2}-1}{q^{4}} $ when only one of the scalers $ \gamma_1 $ or $ \gamma_2 $ are zero. \\
\textbf{Case 2.} $ \dim L^{2}=3. $ \\
The set  $ \lbrace [x_1, x_2], [x_1, x_3], [x_2, x_3]  \rbrace $ is a basis for $ L^{2}. $  
%Therefore $ a=\alpha_1 x_1+\alpha_2 x_2+\alpha_3 x_3+\beta z $ such that  $ z \in Z(L) $ and $ \alpha_i \in \mathbb{F}_{q} $ such that $ 1\leq i \leq 3. $
 Similar to the method of previous case, we have $ \dim \operatorname{Im}ad_{a} =2 $ for all $ a\in L\setminus Z(L) $ and so the number of such elements is equal to $ q^{t}(q^{3}-1). $
Therefore 
\begin{align*}
d(L)&=\dfrac{| Z(L)|}{| L|}+
\dfrac{1}{| L|} \sum_{x\in L\setminus  Z(L)} \frac{1}{| \operatorname{Im}ad_{x}|} \cr
&=\dfrac{q^{t}}{q^{3+t}}+\dfrac{1}{q^{3+t}}(\dfrac{q^{3+t}-q^{t}}{q^{2}})=\dfrac{q^{3}+q^{2}-1}{q^{5}}.
\end{align*}
\end{proof}
The following corollary is a direct consequence of the above theorem and we  omit its proof.
\begin{cor}
Let $ L $ be a Lie algebra over the field $ \mathbb{F}_{2} $ such that $ \dim L/Z(L)=3. $ Then $ d(L) $ is  $ \frac{7}{16} $ or $ \frac{11}{32}. $
\end{cor}
\section{the asymptotic commutativity degree of Lie algebra}
In this section, we define the asymptotic commutativity degree of a family of Lie algebras and prove that there exists  a family of Lie algebras with the asymptotic commutativity degree  $ \frac{1}{q^{k}} $ for all   positive integer $ k. $
\begin{defn}
Let $ \lbrace L_n \rbrace_{n\geq 1} $ be a family of  Lie algebras such that $| L_n|\rightarrow \infty  $ as $ n\rightarrow \infty. $ Then 
$ \lim_{n\to \infty} d(L_n)$ (if exists) is called the asymptotic commutativity degree of $\lbrace L_n\rbrace_{n\geq 1}. $
\end{defn}
Recall from \cite{N2}, 
the Lie algebra $ L $ is  a central product of ideals $ A $ and $ B, $ (denoted by $ A\dotplus B $) if we have $ L=A+B, $ $ [A, B]=0 $ and $ A\cap B\subseteq Z(L). $ This  product is used in the following structures.\\
In \cite[Theorem 4.6]{N2}, the structure of Lie algebras of nilpotency class $ 3 $ with the derived subalgebra of dimension $ 2$   are determined as follows:
\begin{align}\label{struc}
&L_{4,3}=\langle a_1, a_2, a_3, z \mid [a_1, a_2]=a_3, [a_1, a_3]=z \rangle,\cr
& L_{5,5}=\langle a_1, a_2, a_3, a_4, z  \mid [a_1, a_2]=a_3, [a_1, a_3]=z, [a_2, a_4]=z \rangle,\cr
&L_{4,3} \dotplus H(m),\cr
&L_{5,5} \dotplus H(m).
\end{align}
In the following theorem, we determine the commutativity degree of above Lie algebras.
\begin{thm}
Let $ L $ be an $ n $-dimensional Lie algebra of nilpotency class $ 3 $ with the derived subalgebra of dimension $ 2. $ Then
$$
d(L)=
\begin{cases}
\frac{q^{2}-q}{q^{n}}+\frac{q^{2}+q-1}{q^{4}}&\text{if }n \text{~is even, }\\
\frac{q-1}{q^{n}}+\frac{q^{2}+q-1}{q^{4}}&\text{if }n \text{~is odd. }\\
\end{cases}
$$
\end{thm}
\begin{proof}
Since $ L $ is a Lie algebra of nilpotency class $ 3 $ with the derived subalgebra of dimension $ 2,$  $ L $ is isomorphic to one of the Lie algebras 
$ L_{4,3}, $ $ L_{5,5}, $
$ L_{4,3}\dotplus H(m) $ or $ L_{5,5}\dotplus H(m) $ for all  $ m\geq 1 $ by  \eqref{struc}.
By \cite[Example 3.8]{Sh1}, we have $ d(L_{5,5})=\frac{q^{3}+q^{2}-1}{q^{5}} $ and by a similar way in \cite[Example 3.8]{Sh1}, one can see that $ d(L_{4,3})=\frac{2q^{2}-1}{q^{4}}. $
% Theorem $ 2.1 $ in \cite{Sh1}, the commutativity degree of  and $ L_{5,5} $ are equal to $  $ and  $ ,$ respectively. 
Now, let $ L $ be isomorphic to $ L_{5,5}\dotplus H(m) $ for $ m\geq 1. $
Then
\begin{align}\label{00}
d(L)&=\dfrac{1}{| L|}\sum_{x \in L} \dfrac{1}{|  \operatorname{Im}ad_{x}| }=\dfrac{| Z(L)|}{| L|}+ \dfrac{1}{| L|} \sum_{x \in L\setminus Z(L)} \dfrac{1}{|  \operatorname{Im}ad_{x}|}\cr
&=\dfrac{| Z(L)|}{| L|}+ \dfrac{1}{| L|} \sum_{w_1 \in L_{5,5}\setminus Z(L)} \dfrac{1}{|  \operatorname{Im}ad_{x}|}+\dfrac{1}{| L|} \sum_{w_2 \in H(m)\setminus Z(L)} \dfrac{1}{|  \operatorname{Im}ad_{x}|}\cr
&+\dfrac{1}{| L|} \sum_{x=w_1+w_2 \in L\setminus( Z(L)\cup L_{5,5} \cup H(m))} \dfrac{1}{|  \operatorname{Im}ad_{x}|}.\cr
\end{align}
Since  $ \dim L=n $ and $ \dim Z(L)=1, $ then   $\frac{| Z(L)|}{| L|}=\frac{q}{q^{n}}.  $
By the computations of \cite[Example 3.8]{Sh1}, $ \sum_{x \in L_{5,5}\setminus Z(L)} \frac{1}{|  \operatorname{Im}ad_{x}|}=q^{3}+q^{2}-q-1 $ and
$ \sum_{x \in H(m) \setminus Z(L)} \frac{1}{|  \operatorname{Im}ad_{x}|}=q^{2m}-1. $
Also,  $ [L_{5,5}, H(m)]=0$ so,  it is obvious  $   \operatorname{Im}ad_{w_1+w_2}= \langle \operatorname{Im}ad_{w_1}, \operatorname{Im}ad_{w_2}\rangle $ such that $w_1 \in L_{5,5} $ and $ w_2 \in H(m)$ for $m\geq 1. $ 
On the other hand, we can write $ w_1=\alpha_1 a_1+\alpha_2 a_2+ \alpha_3 a_3+\alpha_4 x_4+\alpha_5 z $ such that 
 $ w_1 \in L_{5,5}\setminus Z(L)$  and $ \alpha_i \in \mathbb{F}_q $ for $ 1\leq i \leq 5. $   
 If $ \alpha_1 \neq 0$ or $ \alpha_2 \neq 0, $ then $  \operatorname{Im}ad_{w_1}=\langle x_3, z \rangle $  
 and the number of such elements is $ q^{5}-q^{3}. $
 Otherwise, $  \operatorname{Im}ad_{w_1}=\langle z \rangle$ and the number of them is $ q^{3}-q. $ 
 Moreover, $ \dim H(m)^{2}=1 $ implies that  $  \operatorname{Im}ad_{w_2}=\langle z \rangle $
 for all $ w_2 \in H(m)\setminus Z(L) $ and the number of such elements is equal to $ q^{2m+1}-q. $ Thus the number of elements 
 $ w_1+w_2 \in L\setminus ( Z(L)\cup  L_{4,3}\cup H(m)) $  such that  $ \dim \operatorname{Im}ad_{w_1+w_2}=1 $ or $ 2 $ is equal to 
$(q^{3}-q)(q^{2m}-1) $ or $(q^{5}-q^{3}) (q^{2m}-1), $ respectively.
% Hence the number of  $ w_1+w_2 $ such that $w_1+w_2\in L\setminus \lbrace Z(L), L_{4,3}, H(m)\rbrace $ and  $  \dim \operatorname{Im}ad_{w_1+w_2}=1 $ is equal to $ q^{3}-q $ and 
%the number of  $ w_1+w_2 $ such that $  \dim \operatorname{Im}ad_{w_1+w_2}=2 $ is equal to  $ q^{4}-q^{3}. $  \\ 
% First, let $ L $ is isomorphic to $ L_{4,3}=\langle x_1, x_2, x_3, x_4\mid [x_1, x_2]=x_3, [x_1, x_3]=x_4\rangle. $ By using the proof of Theorem 4.1 in \cite{degree}, we have $ | L_{4,3}| P(L_{4,3})-Z(L_{4,3})=2q^{2}-q-1. $ On the other hand, in the proof of Theorem \ref{} we showed that $ \dim \operatorname{Im}ad_{w_1+w_2}=\dim \operatorname{Im}ad_{w_2} $ such that 
%$ w_1\in H(m)$ for all $ m\geq 0 $ and $w_2 \in L_{4,3} $ 
 %??? Also, $ \dim L_{4,3}\dotplus H(m) $ for all $ m $ such that $ m\geq 0 $ is even.  \\
Therefore, we have $d(L_{5,5}\dotplus (H(m))= \frac{q-1}{q^{n-1}}+\frac{q^{2}+q-1}{q^{4}} $ for  $ m\geq 1, $ by \eqref{00}.
We note that $\dim L=\dim L_{5,5}+ \dim H(m)-1=5+2m $ is odd number.  \\
 If $ L $ is isomorphic to $ L_{4,3}\dotplus H(m) $ for   $ m\geq 1, $ then $ \dim L=\dim L_{4,3}+ \dim H(m)-1=4+2m  $ is even number. By the same computation as above, we can see that  $ d(L)=\frac{q^{2}-q}{q^{n}}+\frac{q^{2}+q-1}{q^{4}} $ in this case. Hence the proof is completed.
\end{proof}
In the following, we give the asymptotic commutativity degree of Lie algebras for the structures  \eqref{struc}.
\begin{cor}
Let $ \lbrace L_n\rbrace_{n\geq 1} $ be a family of  $ n $-dimensional Lie algebras of nilpotency class $ 3 $ with derived subalgebra of dimension $ 2. $ Then 
$| L_n| \rightarrow \infty $ as $ n $ tends to infinity and we have
$$ \lim_{n\to \infty} d(L_n)=\frac{q^{2}+q-1}{q^{4}}.$$
\end{cor}
%In the following theorem, the existence of a family of Lie algebras is shown
% that its asymptotic commutativity degree is  $ \frac{1}{q^{n}} $ for every $ n\geq  1 $ and $ q\geq 2. $ 
In the next theorem, we show that $ \frac{1}{q^{n}} $ is the asymptotic commutativity degree of some family of Lie algebras for every positive integer $ n. $
\begin{thm}
%Let $\lbrace  L_n \rbrace_{n\geq 1}$ be a family of  Lie algebras   with the derived subalgebra of dimension $ 1. $ Then $ \lim_{n\to \infty} d(L_n)=\frac{1}{q}.$ 
%as $ n $ tends to infinity we have $d(L)\rightarrow \frac{1}{q}  $ as $ n\rightarrow \infty. $
There exists a family of Lie algebra having asymptotic commutativity degree $ \frac{1}{q^{n}} $ for every $ n\geq 1. $
\end{thm}
\begin{proof}
Firstly, we show that there is a family of Lie algebras  with the derived subalgebra of dimension $ 1 $ such that the asymptotic commutativity degree $ \frac{1}{q}. $
 Since $ \dim L_n^{2}=1, $ we have  $ d(L_n)= \frac{q^{2n}+q-1}{q^{2n+1}}$ for all $ n\geq 1$ by 
 Corollary \ref{aso1} and we can put $ L_n=H(n)$ by  \cite[Example 3.3]{Sh1}. 
  If $ n $  tends to infinity, then $| L_n|=q^{n} \rightarrow \infty. $ So, $ \lim_{n\to \infty} d(L_n)=\frac{1}{q}.$   Hence by \cite[Proposition 3.5]{Sh1}, we have
 \begin{equation*}
\lim_{n\to \infty} d( \underbrace{L_n\oplus \dots \oplus L_n}_{n_-times})= \lim_{m\to \infty}   (d( L_n))^{n} =\frac{1}{q^{n}}
\end{equation*}
 
  %we have $\lim_{n\to \infty} d( \underbrace{L_n\oplus \dots \oplus L_n}_{n_-times})=\frac{1}{q^{n}}$ 
   when $ n $
  tends to infinity.  
\end{proof}
%\begin{thm}
%There exists a family of Lie algebra having asymptotic commutativity degree $ \frac{1}{q^{n}} $ for every $ n\geq 1. $
%\end{thm}
%\begin{proof}
%Let $ \lbrace L_n \rbrace_{n\geq 1}$ be a family of  $ n $-dimensional Lie algebras over the field $ \mathbb{F}_{q}. $ Then $ | L_n|=q^{n} $ and so $ q^{n}= q^{n_{1}} q^{n_{2}}\dots q^{n_{t}}  $ by the prime factorization of $ n. $ Consider the family of Heisenberg Lie algebras $ H(m) $ for $ m\geq 1. $ By  \cite[Proposition 3.5 and  Example 3.3]{Sh1},
%we have  
%\begin{equation*}
%d( \underbrace{H(m_i)\oplus \dots \oplus H(m_i)}_{n_i-times})= (\frac{1}{q}+\frac{q-1}{q^{2m_i+1}})^{n_i}.
%\end{equation*}
%Now, by \cite[Proposition 3.5]{Sh1} and  considering the following family
%\begin{equation*}
%(H(m_1))^{n_1}\oplus (H(m_2))^{n_2}\dots \oplus  (H(m_t))^{n_t}
%\end{equation*}
%obtained. 
%by Heisenberg algebra. Therefore 
%Since $ \lim_{m_i\to \infty} d((H(m_i))^{n_i})= \frac{1}{q_i^{n_i}} $ for all $ i $
%such that $ 1\leq i\leq t, $  then 
%\begin{align*}
%&\lim_{m_i\to \infty} d( (H(m_1))^{n_1}\oplus (H(m_2))^{n_2}\dots \oplus  (H(m_t))^{n_t})\cr
%&= (\lim_{m_1\to \infty}   d( H(m_1))^{n_1}) \dots (\lim_{m_t \to \infty}   d( H(m_t))^{n_t})=\frac{1}{q^{n}}. 
%\end{align*}
%\end{proof}
It is interesting to see that for rational number $r$  with  $r \in (0, 1) $ there exists at least a family $\lbrace L_n \rbrace_{n\geq 1}$ with the asymptotic commutativity degree. 
%$\lim_{n\to \infty}  d(L_n)=r.$ 
We end the paper with the following open problem.
\begin{theorem}
Let $r$ be an arbitrary rational number in $(0, 1).$ Does there exist a family $\lbrace L_n \rbrace_{n\geq 1}$ of  Lie algebras over the field $\mathbb{F}_q$ such that 
 $ \lim_{n\to \infty} d(L_n)=r$ and $| L_n|\rightarrow \infty  $ as $ n\rightarrow \infty $ ? 
\end{theorem}


\begin{thebibliography}{9}
\bibitem{breath} 
Khuhirun, Borworn. Classification of nilpotent Lie algebras with small breadth. North Carolina State University, 2014.

\bibitem{Les} 
Lescot, Paul. Sur certains groupes finis." Rev. Math. Spéciales 8 (1987): 276-277.


\bibitem{Lescot} 
Lescot, Paul. Degré de commutativité et structure d’un groupe fini. Rev. Math. Spéciales 8 (1988): 276-279. 

\bibitem{Nath} 
Nath, Rajat Kanti. Nath, Rajat Kanti. On asymptotic commutativity degree of finite groups. Proyecciones (Antofagasta) 38.4 (2019): 829-835.

\bibitem{N}
Niroomand, Peyman. Niroomand, Peyman. On dimension of the Schur multiplier of nilpotent Lie algebras. Central European Journal of Mathematics 9.1 (2011): 57-64.

%\bibitem{N1} 
%Johari, Farangis; Parvizi, Mohsen; Niroomand, Peyman. Capability and Schur multiplier of a pair of Lie algebras. J. Geom. Phys. 114 (2017), 184--196.
\bibitem{N2} 
Niroomand, Peyman; Johari, Farangis and Parvizi, Mohsen. Capable Lie algebras with the derived subalgebra of dimension 2 over an arbitrary field. Linear and Multilinear Algebra 67.3 (2019): 542-554
%\bibitem{Mob} 
%N. Jacobson, Basic Algebra I, W. H. Freeman and Co., San Francisco, (1974), 457-465 

\bibitem{Sh1} 
Shamsaki, Afsaneh; Erfanian, Ahmad and Parvizi, Mohsen. On the commutativity degree of a finite-dimensional Lie algebra. Submitted.

\end{thebibliography}
\end{document}